\documentclass[preprint, 3p, times, 12pt]{elsarticle}

\usepackage{amssymb}
\usepackage{amsfonts}
\usepackage{amsthm}
\usepackage{amsmath}
\usepackage{amstext}
\usepackage{mathrsfs}
\usepackage{amscd}
\usepackage{xypic}

\def\mathcal{\mathscr}
\newfont{\aaa}{cmb10 at 19pt}
\newfont{\bbb}{cmb10 at 11pt}

\newtheorem{theorem}{Theorem}[section]
\newtheorem{proposition}[theorem]{Proposition}

\newtheorem{corollary}[theorem]{Corollary}
\newtheorem{definition}[theorem]{Definition}
\newtheorem{example}[theorem]{Example}

\def\v1{\vspace{1mm}}

\def\leq{\leqslant}

\def\geq{\geqslant}

\newcommand{\beq}{\begin{equation}}
\newcommand{\eeq}{\end{equation}}
\newcommand{\bey}{\begin{eqnarray}}
\newcommand{\eey}{\end{eqnarray}}
\newcommand{\beyy}{\begin{eqnarray*}}
\newcommand{\eeyy}{\end{eqnarray*}}

\newcommand{\ra}{\rightarrow}
\newcommand{\dis}{\displaystyle}

\newcommand{\R}{\mathbb{R}}
\newcommand{\Z}{\mathbb Z}
\newcommand{\N}{\mathbb N}
\newcommand{\C}{\mathscr C}

\newcommand{\F}{\mathscr{F}}

\newcommand{\E}{\mathbb E}
\newcommand{\p}{\mathbb P}
\newcommand{\e}{\text{\rm{e}}}

\newcommand{\La}{\Lambda}
\newcommand{\veps}{\varepsilon}

\newcommand{\pb}{\mathscr{P}}

\newcommand{\wt}{\widetilde}

\def\d{\text{\rm{d}}}
\def\E{\mathbb E}
\def\p{\mathbb P}

\def\ll{\mathscr{L}}

\def\S{\mathcal S}

\journal{International Journal of Control}

\begin{document}

\begin{frontmatter}



\title{Viscosity solutions approach to finite-horizon continuous-time Markov decision process}


\author[label1]{Zhong-Wei Liao\footnote{Corresponding author: Zhong-Wei Liao, E-mail: zhwliao@hotmail.com}}

\affiliation[label1]{organization={College of Education for the Future, Beijing Normal University},
            city={Zhuhai},
            postcode={519087},
            country={China}}

\author[label2]{Jinghai Shao}

\affiliation[label2]{organization={Center for Applied Mathematics, Tianjin University},
            city={Tianjin},
            postcode={300072},
            country={China}}

\begin{abstract}
This paper investigates the optimal control problems for the finite-horizon continuous-time Markov decision processes with delay-dependent control policies. We develop compactification methods in decision processes, and show that the existence of optimal policies. Subsequently, through the dynamic programming principle of the delay-dependent control policies, the differential-difference Hamilton-Jacobi-Bellman (HJB) equation in the setting of discrete space was established. Under certain conditions, we give the comparison principle and further prove that the value function is the unique viscosity solution to this HJB equation. Based on this, we show that among the class of delay-dependent control policies, there is an optimal one which is Markovian.
\end{abstract}




\end{frontmatter}


\section{Introduction}
\label{Introduction}

\noindent Continuous-time Markov decision processes (CTMDPs) have been studied intensively due to their rich application in queuing systems, population processes, see, e.g. the monographs \cite{BR11, GS12, GH, PH} and the extensive references therein. From the viewpoint of realistic applications, it is natural to investigate the optimal control problem with delay-dependent controls. The delay caused in the approach of observing the state of the system, making a decision based on this state, and then inputting this decision back into the studied system. However, the system maybe has changed its state at that time. More generally, this control policies are also known as history-dependent control policies, see, for example, \cite{GHH, GHS, GL, H18, KC, PZ11, PL10, Z17}. In this work we develop the viscosity solutions approach of CTMDPs, and it is worth noting that due to the consideration of delay-dependent controls, the controlled system is no longer a Markovian process. For this reason, relevant theoretical tools, such as compactification methods, comparison principle, differential-difference HJB equations and viscosity solutions approach, have also been discussed again.

It is a fundamental problem in the study of MDPs to distinguish the impact on the value function by taking account of all history-dependent policies or of merely Markovian policies. For the discrete-time MDPs in a finite state space, Derman and Strauch \cite[Theorem 2]{DS66} established a basic result which implies that for any history-dependent policy there exists a randomized Markovian policy such that the associated controlled process admits the same marginal state-action distributions. This result also implies that with respect to the criteria of expected discounted, non-discounted costs and expected average costs, the optimization problem over history-dependent policies and over Markovian policies will derive the same value function, see Derman and Strauch \cite{DS66} and Feinberg et al. \cite{FMS}. For the situations of infinite state spaces or unbounded cost functions, more cautious research methods are needed. We constructed an explicit example (see Appendix) to illustrate that if there are no appropriate constraints on the transition probability matrix and cost function, the value function on the history-dependent policies is not equal to that on the Markovian policies. Therefore, discussing appropriate constraints to ensure consistency of the value function across different policy sets is also one of the topics of this article.

As is well known, the expected finite-horizon criterion is a widely used optimality criterion for CTMDPs optimization problems, which has been studied by numerous works, see e.g. \cite{BR11, GS12, GHH, M68, P75, Y78}. For finite-horizon CTMDPs with finite state and action space, Miller \cite{M68} gave a necessary and sufficient condition for the existence of a piecewise constant optimal policy. Subsequently, the state space of CTMDPs had been generalized to denumerable space (cf. \cite{Y78}) and Borel space (cf. \cite{P75}), and the existence of an optimal Markov policy had been proven under the bounded hypothesis of transition rates and cost functions. Recently, Ba\"{u}erle and Rieder \cite{BR11} studies the finite-horizon CTMDPs with Markov polices by a method based on the equivalent transformation from finite-horizon CTMDPs to infinite-horizon discrete-time Markov decision processes. The corresponding optimality equation had been established according to the existing theory on discrete-time Markov decision processes. In addition, Ghosh and Saha \cite{GS12} considered the finite-horizon CTMDPs in Borel state space with bounded transition rates and Markov policies. The existence of a unique solution to the optimality equation is guaranteed by the Banach fixed point theorem, relatively, the existence of an optimal Markov policy is based on It\^{o}-Dynkin's formula. The finite-horizon CTMDPs with unbounded transition rates are investigated in Guo et al. \cite{GHH}.

The work \cite{GHH} also studied the history-dependent control problem for jumping processes. The precise construction of such kind of controlled system is presented. However, via the main result \cite[Theorem 4.1]{GHH}, the value function $V^\ast (t,i)$ ($t>0$), defined in \cite[p.1069]{GHH},  associated with the optimal control problem over the set of randomized Markov policies can be characterized as a unique solution to a differential equation, and in such case the optimal Markov control policies are shown to exist. Nevertheless, if considering the control problem over the set of history-dependent control policies, there is no result in \cite{GHH} on the existence of the optimal control and on the characterization of the associated value function. In the current work,  we shall show the existence of the optimal controls over the set  of delay-dependent controls and characterize the associated value function.

The approaches used in the aforementioned works in the study of CTMDPs rely on the characterization of the Markov chains, and are not suitable to our current situation any longer since the controlled process is no longer a Markovian one caused by the delays. We develop the compactificion method used usually in the control problem for diffusion processes to the setting of jumping processes in order to show the existence of the optimal delay-dependent control policies. This is the starting point of this work. Precisely, the main contributions of the present paper are as follows:

(i) In comparison with \cite{GS12,GHH}, our method used in the existence of an optimal delay-dependent control does not involve the solvability of the optimality equation, but is based on the compactification method, which is an effective method in the research of the optimal control problem of jump-diffusion processes, cf. \cite{CMR,DM,HS95,HS95b}. The basic idea is inspired by Kushner \cite{Ku}, Haussmann and Suo \cite{HS95,HS95b}. Our approach is also suitable to other optimality criteria in the study of CTMDPs such as expected discounted, average and risk-sensitive.

(ii) According to the measurable selection theorem (cf. Stroock and Varadhan \cite{SV79}), the dynamic programming principle is established in Theorem \ref{t2}, which deduces that the value function is a solution to a HJB equation provided the value function to be regular enough. Here the HJB equation is a differential-difference equation. We develop the viscosity solution approach to such equation, and especially we establish the comparison principle for such differential-difference HJB equation.

The rest of our paper is organized as follows. In Section 2, we state the concept of delay-dependent controls and the optimality problems of CTMDPs, and further introduced the main assumptions of this article. For the convenience, the optimality problem is reformulated on the canonical path space. In Section 3, by developing the compactification method within the framework of MDPs, we prove the existence of the optimal delay-dependent controls. In Section 4, we study the HJB equation derived from the dynamic programming principle through the viscosity solution approach. In order to prove the existence and uniqueness of viscosity solution, we also prove the comparison principle under the framework of MDPs. Invoking the corresponding results on the optimal control problem over Markovian control policies, we further show that there must exist an optimal Markovian control policy for the control problem over the class of delay-dependent control policies.

\section{Formulation and assumptions}
\label{Formulation and assumptions}

\noindent The objective of this section is to describe briefly the controlled process and the associated optimal control criterion in this paper. Let $(\Omega, \F, (\F_t)_{t \geq 0}, \p)$ be a filtered probability space satisfying the usual conditions, i.e. $(\Omega, \F, \p)$ is complete, the filtration $(\F_t)_{t \geq 0}$ is right-continuous and $\F_0$ contains all $\p$-null sets in $\F$. Let $\S=\{1,2,\ldots\}$ be the countable state space, $U$ be the action space which is a compact subset of $\R^k$ for some $k\in\N$. Denote by $\pb(U)$ the collection of all probability measures over $U$, which is endowed with $L_1$-Wasserstein distance $W_1$ defined by:
\begin{equation*}
  W_1(\mu,\nu)=\inf\Big\{\int_{U\times U}\!\!|x-y|\pi(\d x,\d y);\ \pi\in \C(\mu,\nu)\Big\},
\end{equation*}
where $\C (\mu,\nu)$ stands for the set of all couplings of $\mu$ and $\nu$ in $\pb(U)$.
Since $U$ is compact, $\pb(U)$ becomes a compact Polish space under the metric $W_1$, and the weak convergence of probability measures in $\pb(U)$ is equivalent to the convergence in the $W_1$ distance (cf. e.g. \cite[Chapter 7]{AGS}). In this work we investigate the finite-horizon optimal control problem on $[0, T]$, where $T > 0$ is fixed throughout this work.

For each $\mu\in \pb(U)$, $(q_{ij}(\mu))$ is a transition rate matrix over the state space $\S$, which is assumed to be conservative, i.e.
\[\sum_{j\neq i} q_{ij}(\mu)=q_i(\mu)=-q_{ii}(\mu),\quad\forall\,i\in\S,\ \mu\in \pb(U).\]
The process $(\La_t)$ is an $\F_t$-adapted jump process on $\S$ satisfying
\begin{equation}\label{b1}
\p(\La_{t+\delta}=j|\La_{t}=i,\mu_t=\mu)=\begin{cases}
                                           q_{ij}(\mu)\delta+o(\delta), & \mbox{if $i\neq j$ }, \\
                                           1+q_{ii}(\mu)\delta+o(\delta), & \mbox{otherwise},
                                         \end{cases}
\end{equation}
provided $\delta>0$.

In order to introduce the delay-dependent control, we first introduce some notations. Given any metric space $E$, denote by $\mathcal{C}([0,T];E)$ the collection of continuous functions $x: [0,T] \ra E$, and $\mathcal{D}([0,T];E)$ the collection of right-continuous functions with left limits $\lambda:[0,T]\ra E$.  For $r_0\in (0,T)$ and $s\in [0,T]$, define a shift operator $\theta_{s,r_0} : \mathcal{D}([0,T];\S)\ra \mathcal{D}([0,T];\S)$ by
\begin{equation}\label{b2}
(\theta_{s,r_0} \lambda) (t)=\lambda ((t-r_0)\vee s ), \quad t\in [0,T].
\end{equation}
Moreover, $\dis \theta_{s,r_0}^k\lambda (t):=\lambda ((t-kr_0)\vee s)$ for $\lambda\in \mathcal{D}([0,T];\S)$, $k \in \mathbb{Z}_+$. Next, we introduce the concept of delay-dependent control.

\begin{definition}\label{def-1}
Fix an arbitrary $m \in \mathbb{Z}_+$ and $r_0 > 0$. Given any $s\in [0,T)$ and $i\in \S$, a randomized delay-dependent control is a term  $\alpha=(\La_t,\mu_t, s, i)$ such that
  \begin{itemize}
    \item[$(i)$] $(\La_t)$ is an $\F_t$-adapted jump process satisfying \eqref{b1} with initial value $\La_s=i$.
    \item[$(ii)$] There exists a measurable map $h:[0,T]\times \S^{m+1}\to \pb(U)$ such that
    \begin{equation}\label{mu}
    \mu_t=h(t,\theta_{s,r_0}^0\La(t),\ldots, \theta_{s,r_0}^m\La(t))  \quad \text{for almost all}\  t\in [s,T].
    \end{equation}
  \end{itemize}
\end{definition}

The parameter $r_0 > 0$ is used to characterize the time interval of delay of the controlled processes, and $m \in \mathbb{Z}_+$ for the number of delay. The collection of all delay-dependent control  $\alpha$ with initial condition $(s,i)$ is denoted by $\Pi_{s,i}$. When the starting time of the optimal control problem is $s$, as we have no further information on the controlled system before the initial time $s$, we use the state of the process $(\La_t)$ at time $s$ to represent its states before time $s$, which is reflected by the definition of $\mu_t$ through equation \eqref{mu}. Such treatment has been used in the study of optimal control problem over history-dependent policies; see, for instance, \cite{GHH,GHS}.

Let $f:[0,T]\times \S\times \pb(U)\ra [0,\infty)$, $g:\S\ra [0,\infty)$ be two lower semi-continuous functions. The expected cost for the delay-dependent control  $\alpha\in \Pi_{s,i}$ is defined by
\begin{equation}\label{b5}
J(s,i,\alpha)=\E\Big[\int_s^Tf(t,\La_t,\mu_t)\d t+g(\La_T)\Big],
\end{equation}
and the value function is defined by
\begin{equation}\label{b6}
  V(s,i)=\inf_{\alpha\in \Pi_{s,i}} J(s,i,\alpha).
\end{equation}
It immediately implies that the value function $V$ satisfies $V(T, i) = g(i)$, $\forall i \in \S$. A delay-dependent control  $\alpha^\ast\in \Pi_{s,i}$ is said to be \emph{optimal}, if $V(s,i)=J(s,i,\alpha^\ast)$.

The set of delay-dependent controls  introduced in Definition \ref{def-1} contains many interesting control policies. Next, we present some examples below.
\begin{example}\label{ex-1}
We consider the optimal control problem with initial time $s=0$.
\begin{enumerate}
  \item $\mu_t=h(\La_t)$ for some $h:\S\ra \pb(U)$. In this situation, $\alpha$ is corresponding to the stationary randomized Markov policy studied by many works; see, e.g. \cite{GH}.
  \item $\mu_t=h(\La_{(t-r_0)\vee 0})$ for some $h:\S\ra \pb(U)$. Now the control policies are purely determined by the jump process with a positive delay. This kind of controls is very natural to be used in the realistic application.
  \item $\mu_t=h(t,\La_{(t-r_0)\vee 0}, \La_{(t-2r_0)\vee 0})$ for some   $h:[0,T]\times \S\times \S\ra \pb(U)$.
  \item $\mu_t=h(t,\La_{(t-r_0)\vee 0})$ for some $h(t,i)=\delta_{u_t(i)}$ for each $i\in \S$, where $t\mapsto u_t(i)$ is a curve in $U$ and $\delta_x$ denote the Dirac measure in $U$. 
\end{enumerate}
\end{example}

In this paper we impose the following assumptions on the primitive $Q$-matrix of the continuous-time Markov decision process $(\La_t)$.

\noindent\textbf{Assumptions}:
\begin{itemize}
  \item[(H1)] $\mu\mapsto q_{ij}(\mu)$ is continuous for every $i,\,j\in\S$, and
  $$
  M:=\sup_{i\in\S} \sup_{\mu\in \pb(U)} q_i(\mu)<\infty.
  $$
  \item[(H2)] There exists a compact function $\Phi:\S\ra [1,\infty)$, a compact set $B_0\subset \S$, constants $\lambda_0>0$ and $\kappa_0\geq 0$ such that
      \begin{align*}Q_\mu\Phi(i)&:=\sum_{j\neq i}q_{ij}(\mu)\big(\Phi(j)-\Phi(i)\big)\leq \lambda_0\Phi(i)+\kappa_0\mathbf{1}_{B_0}(i).
      \end{align*}
  \item[(H3)] There exists $K\in \N$ such that for every $i\in\S$ and $\mu\in \pb(U)$, $q_{ij}(\mu)=0$, if $|j-i|>K$.
\end{itemize}

Here if for every $c\in \R$, the set $\{i\in\S; \Phi(i)\leq c\}$ is a compact set, then  $\Phi$ is called a compact function. Condition (H3) is a technical condition, which is used when we consider to use the dominated convergence theorem in the argument of our main theorem.

In contrast to the well-studied continuous-time Markov decision process, the controlled system $(\La_t)$ studied in this work is no longer a Markov chain, and the delay-dependent control  policy makes it more difficult to describe the evolution of $(\La_t)$. Following \cite{Sh19b}, we shall develop the classical compactness method to deal with the control problem with delay-dependent controls. The compactification method is usually used to cope with the optimal control problem for stochastic differential equations (cf. Kushner \cite{Ku}, Haussmann and Suo \cite{HS95,HS95b} and references therein). We extend this method to deal with stochastic processes in discrete space.

Let
\begin{equation}\label{a-1}
\mathscr{U}=\{\mu:[0,T]\to \pb(U) \ \text{is measurable}\}.
\end{equation}
$\mathscr{U}$ can be viewed as a subspace of $\pb([0,T]\times U)$ through the map
\[(\mu_t)_{t\in[0,T]}\mapsto \bar\mu,\]
where $\bar\mu$ is determined by
\[\bar\mu(A\times B)=\frac1T\int_A\mu_t(B)\d t.\]
Endow $\mathscr{U}$ with the induced weak convergence topology from $\pb([0,T]\times U)$.
This topology is equivalent to the topology induced by the following Wasserstein distance on $\pb([0,T]\times U)$:
\[W_1(\bar\mu,\bar \nu)=\inf_{\Gamma\in \C(\bar\mu,\bar\nu)}\!\int_{([0,T]\,\times\,U)^2}\!\!\big(|s-t|+|x-y|\big)\d \Gamma((s,x),(t,y)),\]
where $\C(\bar\mu,\bar\nu)$ stands for the collection of couplings of $\bar\mu$ and $\bar\nu$ over $([0,T]\times U)^2$.
The canonical path space for our problem is defined as
$$
\hat{\Omega} = \mathcal{D}([0,T]; \S) \times \mathscr{U}
$$
endowed with the product topology, which is a metrizable and separable space (cf. \cite{HS95}). Denote by $\tilde{\mathcal{D}}^1$ (resp. $\tilde{\mathcal{D}}^2$) the Borel $\sigma$-algebra of $\mathcal{D}([0,T]; \S)$ (resp. $\mathscr{U}$), and $\tilde{\mathcal{D}}^1_t$ (resp. $\tilde{\mathcal{D}}^2_t$) the $\sigma$-algebra up to time $t$. Define the $\sigma$-algebra of $\hat{\Omega}$ as
$$
\hat{\F} := \tilde{\mathcal{D}}^1 \times \tilde{\mathcal{D}}^2, \quad \text{and} \quad \hat{\F}_t = \tilde{\mathcal{D}}^1_t \times \tilde{\mathcal{D}}^2_t.
$$
For each delay control $\alpha = (\Lambda_t, \mu_t,x,i) \in \Pi_{s, i}$, we define a measurable map $\Phi_{\alpha} : \Omega \to \hat{\Omega}$ as
$$
\Phi_{\alpha} (\omega) = (\Lambda_t (\omega), \mu_t (\omega))_{t \in [0, T]}, \quad \Lambda_r (\omega) \equiv i, \ \mu_r (\omega) \equiv \mu_s, \ 0\leq r \leq s.
$$
Then, there exists a corresponding probability on $(\hat{\Omega}, \hat{\F})$ defined by $R= \mathbb{P} \circ \Phi_{\alpha}^{-1}$. We denote by $\hat{\Pi}_{s, i}$ the space of probabilities induced by the delay-dependent control  set $\Pi_{s, i}$ with initial condition $(s, i)$. By the definition of value function, we have
$$
V(s, i) = \inf_{\alpha \in \Pi_{s, i}} J (s, i, \alpha) = \inf_{R \in \hat{\Pi}_{s, i}} \mathbb{E}_{R} \left[ \int_s^T f (t, \Lambda_t, \mu_t) \d t + g (\Lambda_T) \right].
$$
The topology and properties of the canonical path space have been well studied, see, for instance \cite{HS95,MZ84,SV79} and the references therein.

\section{Existence of optimal delay-dependent controls}
\label{Existence of optimal delay-dependent controls}

\noindent By developing the compactification method presented, for instance, in \cite{HS95} and \cite{Ku}, Shao \cite{Sh19b} investigated the optimal control problem for the regime-switching processes. There, the control on the transition rate matrix of the jumping process $(\La_t)$ has been studied. In this paper we shall apply the result \cite[Theorem 2.3]{Sh19b} to the current situation to obtain the existence of optimal delay-dependent controls  of our continuous-time Markov decision processes under the mild conditions $\mathrm{(H1)}$-$\mathrm{(H3)}$.

\begin{theorem}\label{t1}
Assume $\mathrm{(H1)}$-$\mathrm{(H3)}$ hold. Then for every $s\in [0,T)$, $i\in\S$, there exists an optimal delay-dependent control  $\alpha^\ast\in \Pi_{s,i}$.
\end{theorem}

\begin{proof} This theorem is proved by using the idea of \cite[Theorem 2.3]{Sh19b}. The proof is a little long. In order to save space, here we only sketch the idea and point out the different points compared with that of \cite[Theorem 2.3]{Sh19b}.

We only need to consider the nontrivial case $V(s,i)<\infty$. For simplicity of notation, we consider the case $s=0$, and separate the proof into three steps.

\textbf{Step 1.} According to the definition of $V(0,i)$, there exists a sequence of delay-dependent controls  $\alpha_n =(\La_t^{(n)},\mu_t^{(n)}, 0, i)\in \Pi_{0,i}$ such that
\begin{equation}\label{c-1}
\lim_{n\ra \infty} J(0,i,\alpha_n)=V(0,i).
\end{equation}

Denote by $R_n$ the probability measures on $(\hat{\Omega}, \hat{\F})$ corresponding to $\alpha_n$. Let $\ll_\mu^n$ (resp. $\ll_\La^n$) be the marginal distribution of $R_n$ with respect to $(\mu_t^{(n)})_{t\in [0,T]}$ (resp. $(\La_t^{(n)})_{t\in [0,T]}$) in $\mathscr{U}$ (resp. $\mathcal D([0,T]; \S)$). Since $\pb([0,T]\times U)$ is compact and further $\mathscr{U}$ is compact as a closed subset, we have $(\ll_\mu^n)_{n\geq 1} $ is tight.

We proceed to prove that $(\ll_\La^n)_{n\geq 1}$ is tight. For each $n \geq 1$, by (H2) and It\^o-Dynkin's formula (cf. \cite[Theorem 3.1]{GHH}), we have
\begin{align*}
  \E \Phi(\La_t^{(n)})&=\Phi(i)+\E\int_0^t Q_{\mu_s}\Phi(\La_s^{(n)})\d s\\
         &\leq \Phi(i)+\E\int_0^t \big(\lambda_0\Phi(\La_s^{(n)})+\kappa_0\big)\d s,
\end{align*}
which yields from Gronwall's inequality that
\begin{equation}\label{c-8}
\E \Phi(\La_t^{(n)})\leq \big(\Phi(i)+\kappa_0 T\big)\e^{\lambda_0 t},\quad t\in [0,T], \  n \geq 1.
\end{equation}
For any $\veps>0$, we can find $N_\veps>0$ such that
\begin{equation}\label{c-9}
\sup_n\p(\La_t^{(n)}\in K_\veps^c)\leq \sup_{n}\frac{\E \Phi(\La_t^{(n)})}{N_\veps}\leq \frac{(\Phi(i)+\kappa_0 T)\e^{\lambda_0 T}}{N_\veps}<\veps,
\end{equation}
where $K_\veps=\{j\in \S; \Phi(j)\leq N_\veps\}$. Since $\Phi$ is a compact function, $K_\veps$ is a compact set. Moreover, for every $0\leq u\leq \delta$, due to (H1),
\begin{equation}\label{c-10}
\begin{split}
\E\big[\mathbf{1}_{\La_{t+u}^{(n)}\neq \La_{t}^{(n)}}\big]&\leq 1-\p(\La_s^{(n)}=\La_t,\forall\, s\in [t, t+u])\\
&\leq 1-\e^{-Mu}\leq 1-\e^{-M\delta}=:\gamma_n(\delta).
\end{split}
\end{equation}
To apply \cite[Theorem 8.6, p.138]{EK}, by taking $q(i,j)= \mathbf{1}_{i\neq j}$, $\beta=1$, $\gamma_n(\delta)$ given in \eqref{c-10} and invoking \eqref{c-9}, we obtain the tightness of $(\ll_\La^n)_{n\geq 1}$.

\textbf{Step 2.} Since the marginal distributions $(\ll_\La^n)_{n\geq 1}$ and $(\ll_\mu^n)_{n\geq 1}$ are both tight, $(R_n)_{n \geq 1}$ is tight as well. Hence, there exists a subsequence $n_k$, $k\geq 1$, such that $R_{n_k}$ weakly converges to some probability measure $R_0$ on $(\hat{\Omega}, \hat{\F})$ as $k\ra \infty$. By virtue of Skorokhod's representation theorem (cf. e.g. \cite[Chapter 3]{EK}), there exists a probability space $(\Omega',\F',\p')$ on which is defined a sequence of $\hat{\Omega}$-valued random variables $Y_{n_k}=(\La_t^{(n_k)},\mu_t^{(n_k)})_{t\in [0,T]}$ with distribution $R_{n_k}$, $k \geq 1$, and $Y_0=(\La_t^{(0)},\mu_t^{(0)})_{t\in[0,T]}$ with distribution $R_0$ such that
\begin{equation}\label{c-11}
\lim_{k\ra \infty} Y_{n_k}=Y_0,\quad \p'\text{-a.s.}
\end{equation}

Analogous to the Step 2 in the argument of \cite[Theorem 2.3]{Sh19b}, we can show that $\alpha^* := (\La_t^{(0)},\mu_t^{(0)}, 0, i)$ is a delay-dependent control in $\Pi_{0,i}$. During this procedure, we need to replace the sigma fields $\mathscr{F}^{X,\La}_{-n, t}$ by the following
\[\mathscr{F}^{\La}_{-n,t}:=\overline{\sigma\{(\La_t^{(k)},\ldots,   \La_{t-mr_0}^{(k)}); k\geq n\}}. \]

\textbf{Step 3.}
Invoking \eqref{c-1} and the lower semi-continuity of $f$ and $g$, we obtain
\begin{align*}
  V(0,i)&=\lim_{k\ra \infty} \E\Big[\int_0^T f(t,\La_t^{(n_k)},\mu_t^{(n_k)})\d t+g(\La_T^{(n_k)}) \Big]\\
  & \geq \E\Big[\int_0^T f(t, \La_t^{(0)},\mu_t^{(0)})\d t+g(\La_T^{(0)})\Big]\\
  &\geq V(0,i).
\end{align*}
By taking $\alpha^\ast=(\La_t^{(0)},\mu_t^{(0)}, 0, i) \in \Pi_{0, i}$, the previous inequalities imply that $\alpha^\ast$ is an optimal delay-dependent control  of the continuous-time Markov jump process. The proof of this theorem is complete.
\end{proof}

\section{Dynamic programming principle and viscosity solution}
\label{Dynamic programming principle and viscosity solution}

\noindent In the rest of the paper, we introduce the dynamic programming principle for the controlled processes with delay-dependent control  and the differential equation satisfied by the value function. To do so, we introduce some notations. Assume that $\tau$ is an $\hat{\F}_t$-stopping time satisfying $0 \leq \tau \leq T$, $\hat{\F}_{\tau}$ is denoted by the collection of sets $A$ such that $A \cap \{ \tau \leq t\} \in \hat{\F}_t$, $\forall t \in [0, T]$.

\begin{theorem}\label{t2}
Assume $\mathrm{(H1)}$-$\mathrm{(H3)}$  hold. For each $\hat{\F}_t$-stopping time $\tau$ satisfying $s \leq \tau \leq T$, then
$$
V(s, i) = \inf \left\{ \mathbb{E}_{R} \left[ \int_s^\tau f (t, \Lambda_t, \mu_t) \d t + V(\tau, \Lambda_{\tau}) \right] ; R \in \hat{\Pi}_{s, i} \right\}.
$$
\end{theorem}

\begin{proof} Define a subset of $\hat{\Pi}_{s, i}$ as
$$
\hat{\Pi}^0_{s, i} = \left\{ R \in \hat{\Pi}_{s, i}: V(s, i) = \mathbb{E}_{R} \left[ \int_s^T f (t, \Lambda_t, \mu_t) \d t + g (\Lambda_T) \right] \right\}.
$$
By Theorem \ref{t1}, $\hat{\Pi}^0_{s, i} \neq \emptyset$ for any $s \in [0, T]$ and $i \in \S$. According to measurable choices theorem presented by Stroock and Varadhan \cite{SV79}, there exists a Borel-measurable map $H: [0, t] \times \S \to \pb(U)$, which is called measurable selector, satisfying for each $(s, i) \in [0, t] \times \S$, $H(s, i) \in \hat{\Pi}^0_{s, i}$. Refer to \cite[Lemma 3.9]{HS95} for more details of the existence of the measurable selector. Hence, for any $\hat{\omega} \in \hat{\Omega}$, $H (\tau(\hat{\omega}), \Lambda_{\tau(\hat{\omega})})$ is a probability measure  on $(\hat{\Omega}, \hat{\F})$ and satisfies
\begin{equation}\label{d-2}
V(\tau(\hat{\omega}), \Lambda_{\tau(\hat{\omega})}) = \mathbb{E}_{H (\tau(\hat{\omega}), \Lambda_{\tau(\hat{\omega})})} \left[ \int_{\tau(\hat{\omega})}^T f (t, \Lambda_t, \mu_t) \d t + g(\Lambda_T) \right].
\end{equation}
Note that the topology on $\hat{\Omega}$ is separable, then $\hat{\F}_t$ is countably generated, and then for every probability measure $\p$ on $(\hat{\Omega}, \hat{\F})$, the regular conditional probability distribution of $\p$ for given $\hat{\F}_\tau$ exists, cf. \cite{HS95,HS95b}. According to \cite[Lemma 3.3]{HS95b}, for each $R \in \hat{\Pi}_{s, i}$, there exists a unique probability measure, denoted by $R^H$, such that $R^{H} (A) = R (A)$, $\forall \ A \in \F_{\tau}$ and the regular conditional probability distribution of $R^H$ for given $\F_{\tau}$ is $H (\tau(\cdot), \Lambda_{\tau(\cdot)})$. Moreover, by \cite[Proposition 3.8]{HS95b}, it holds that $R^{H} \in \hat{\Pi}_{s, i}$. Hence, we have
$$
V(\tau(\hat{\omega}), \Lambda_{\tau(\hat{\omega})}) = \mathbb{E}_{R^H} \left[ \int_{\tau(\hat{\omega})}^T f (t, \Lambda_t, \mu_t) \d t + g(\Lambda_T) \Big| \F_{\tau} \right].
$$

Due to the definition of value function $V(s, i)$, we have
\begin{align*}
V(s, i) & \leq \mathbb{E}_{R^H} \left[ \int_s^\tau f (t, \Lambda_t, \mu_t) \d t + \int_\tau^T f (t, \Lambda_t, \mu_t) \d t + g(\Lambda_T) \right] \\
&= \mathbb{E}_{R^H} \left[ \int_s^\tau f (t, \Lambda_t, \mu_t) \d t + \mathbb{E}_{R^H} \left[ \int_\tau^T f (t, \Lambda_t, \mu_t) \d t + g(\Lambda_T) \Big| \F_{\tau} \right] \right] \\
&= \mathbb{E}_{R} \left[ \int_s^\tau f (t, \Lambda_t, \mu_t) \d t + V(\tau, \Lambda_{\tau}) \right],
\end{align*}
where the last equation is based on the relationship between $R$ and $R^H$. The arbitrariness of $R \in \hat{\Pi}_{s, i}$ implies that
$$
V(s, i) \leq \inf \left\{ \mathbb{E}_{R} \left[ \int_s^\tau f (t, \Lambda_t, \mu_t) \d t + V(\tau, \Lambda_{\tau}) \right] ; R \in \hat{\Pi}_{s, i} \right\}.
$$

Conversely, by Theorem \ref{t1}, there exists an optimal delay-dependent control  $\alpha^* \in \Pi_{s, i}$ and then denote by $R^* \in \hat{\Pi}_{s, i}$ the corresponding probability measure on $(\hat{\Omega}, \hat{\F})$. Then we have
\begin{align*}
V(s, i) & = \mathbb{E}_{R^*} \left[ \int_s^\tau f (t, \Lambda_t, \mu_t) \d t + \int_\tau^T f (t, \Lambda_t, \mu_t) \d t + g(\Lambda_T) \right] \\
&\geq \mathbb{E}_{R^*} \left[ \int_s^\tau f (t, \Lambda_t, \mu_t) \d t + V(\tau, \Lambda_{\tau}) \right] \\
&\geq \inf \left\{ \mathbb{E}_{R} \left[ \int_s^\tau f (t, \Lambda_t, \mu_t) \d t + V(\tau, \Lambda_{\tau}) \right] ; R \in \hat{\Pi}_{s, i} \right\}.
\end{align*}
The dynamic programming principle is thus proved.
\end{proof}

The next result is about the continuity of value function. Since $\S$ is a countable state space equipped with discrete topology, we only need to consider the continuity of $V(s, i)$ in the time variable $s$.

\begin{proposition}\label{prop1}
Assume $\mathrm{(H1)}$-$\mathrm{(H3)}$ hold. Suppose that $f$, $g$ are bounded and $f$ satisfies the following condition,
\begin{equation}\label{d-0}
|f(t, i, \mu) - f(s, i, \mu)| \leq C_0 |t- s|, \qquad 0 \leq s, t \leq T,
\end{equation}
uniformly for $i \in \S$ and $\mu \in \pb(U)$. Then, the value function $V(s, i)$ is Lipschitz continuous with respect to the time variable $s$. In fact, there exists a constant $C > 0$ such that for any $i \in \S$
$$
\left| V(s, i) - V(s', i) \right| \leq C |s -s'|, \quad 0 \leq s,  s' \leq T.
$$
\end{proposition}

\begin{proof} For convenience, denote by $C_1$ and $C_2$ the constants such that
$$
\sup_{(t, i, \mu) \in [0, T] \times \S \times \pb(U)} |f (t, i, \mu)| \leq C_1 \quad \text{and} \quad \sup_{i \in \S} |g (i)| \leq C_2.
$$
Fix any $i \in \S$ and assume $0 \leq s \leq s' \leq T$. According to Theorem \ref{t1}, there exists an optimal delay-dependent control  $\alpha^* = (\Lambda_t, \mu_t, s, i) \in \Pi_{s, i}$ such that $V(s, i) = J(s, i, \alpha^*)$. By time shift, we can define a couple of processes with initial point $(s', i)$ as following
$$
\Lambda'_t = \Lambda_{t - \Delta s}, \quad \mu'_t = \mu_{t- \Delta s}, \quad \forall\, t \in [s', T],
$$
where $\Delta s := s' - s$. It is easy to verify that (\ref{b1}) and (\ref{mu}) hold for $(\Lambda'_t, \mu'_t)$, which means that $\alpha' := (\Lambda'_t, \mu'_t, s', i)$ is a delay-dependent control  in $\Pi_{s', i}$. Using (H1) and (\ref{b1}), we have
$$
\mathbb{E} \left[ \mathbf{1}_{\Lambda'_t \neq \Lambda_t} \right] = \mathbb{P} \left( \Lambda_{t- \Delta s} \neq \Lambda_{t} \right) \leq M \Delta s + o (\Delta s).
$$
By the definition of the value function, we have
\begin{align}\label{d-1}
|V (s', i) - V(s, i) | & \leq \mathbb{E} \left[ \int_{s'}^T \left| f(t, \Lambda'_t,  \mu'_t) - f (t, \Lambda_t, \mu_t) \right| \d t\right] \notag \\
& \quad + \mathbb{E} \left[ \left| g(\Lambda'_{T}) - g (\Lambda_T) \right| \right] + \mathbb{E} \left[\int_{s}^{s'} \left| f (t, \Lambda_t, \mu_t) \right| \d t\right].
\end{align}
According to the boundedness of $f$ and $g$, we obtain
\begin{align*}
\mathbb{E} & \left[ \int_s^{s'} |f (t, \Lambda_t, \mu_t)| \d t \right] \leq C_1 \Delta s, \quad \text{and}\\
\mathbb{E} & \left[ |g (\Lambda'_T ) - g (\Lambda_T)| \right] \leq 2 C_2 \mathbb{E} \left[\mathbf{1}_{\Lambda'_T \neq \Lambda_T}\right] \leq 2 M C_2 \Delta s + o (\Delta s).
\end{align*}
To estimate the first term of (\ref{d-1}), we combine the boundedness and (\ref{d-0}),
\begin{align*}
& \mathbb{E} \left[ \int_{s'}^T | f (t, \Lambda'_t, \mu'_t) - f (t, \Lambda_t, \mu_t) | \d t \right] \\
& = \mathbb{E} \left[ \int_s^{T - \Delta s } | f (t+ \Delta s, \Lambda_t, \mu_t) \d t - f (t, \Lambda_t, \mu_t) | \d t \right] \\
& \quad + \mathbb{E} \left[ \int_s^{s'} | f(t, \Lambda_t, \mu_t) | \d t \right] + \mathbb{E} \left[ \int_{T- \Delta s}^T | f (t, \Lambda_t, \mu_t) | \d t \right] \\
& \leq T C_0 \Delta s + 2 C_1 \Delta s.
\end{align*}
Hence,
$$
|V(s, i) - V(s', i)| \leq (3 C_1 + 2 M C_2 + T C_0) \Delta s + o (\Delta s).
$$
By the symmetric position of $s$  and $s'$, we have $|V(s, i) - V(s', i)| \leq C |s - s'|$.
\end{proof}

According to Proposition 4.2 and Rademacher's theorem, we know that for each $i \in \S$, $t \to V(t, i)$ is almost everywhere differentiable in $[0, T]$ with respect to Lebesgue measure. In some practical applications, the property of almost everywhere differentiable is not enough, especially when $\S$ is a general state space rather than a countable space. But, it is not easy to justify whether $V (t,i)$ is differentiable every where in $[0,T]$. In such situation, it is useful to introduce the concept of viscosity solution to further characterize $V (t,i)$. Consider the following equation
\begin{equation}\label{HJB0}
-\frac{\partial v}{\partial t}-\inf_{\mu\in \pb(U)} \Big\{\sum_{j\neq i}q_{ij}(\mu)\big(v(t,j)-v(t,i)\big)+f(t,i,\mu)\Big\}=0,\quad v(T,i)=g(i).
\end{equation}

\begin{definition}\label{def-2}
Let $v:[0,T)\times \S\to \R$ be a continuous function.
\begin{itemize}
  \item[$(i)$] $v$ is called a viscosity supersolution of \eqref{HJB0} if $\phi(T,i)=g(i)$,
  \[-\frac{\partial \phi}{\partial t}(t_0,i_0)-\inf_{\mu\in \pb(U)} \Big\{\sum_{j\neq i_0}q_{i_0j}(\mu)\big(\phi(t_0,j)-\phi(t_0,i_0)\big)
  +f(t_0,i_0,\mu)\Big\}\geq 0\]
  for all $(t_0,i_0)\in [0,T)\times \S$ and for all $\phi\in C^1([0,T)\times \S)$ such that $(t_0,i_0)$ is a minimum point of $v-\phi$.
  \item[$(ii)$] $v$ is called a viscosity subsolution of \eqref{HJB0} if $\phi(T,i)=g(i)$,
  \[-\frac{\partial \phi}{\partial t}(t_0,i_0)-\inf_{\mu\in \pb(U)} \Big\{\sum_{j\neq i_0}q_{i_0j}(\mu)\big(\phi(t_0,j)-\phi(t_0,i_0)\big)
  +f(t_0,i_0,\mu)\Big\}\leq 0\]
  for all $(t_0,i_0)\in [0,T)\times \S$ and for all $\phi\in C^1([0,T)\times \S)$ such that $(t_0,i_0)$ is a maximum point of $v-\phi$.
  \item[$(iii)$] $v$ is called a viscosity solution to \eqref{HJB0} if it is both a viscosity subsolution and a viscosity supersolution of \eqref{HJB0}.
\end{itemize}
\end{definition}

The next result says that the value function is a solution to the HJB equation \eqref{HJB0} in the viscosity sense.

\begin{theorem}\label{t3}
Under the conditions of Proposition \ref{prop1}, the value function $V(t,i)$ is a viscosity solution to   \eqref{HJB0}.
\end{theorem}

\begin{proof} We first consider the viscosity subsolution property. Let $(t_0,i_0)\in [0,T)\times \S$ and $\phi\in C^1([0,T)\times \S)$ be a test function such that
\begin{equation}\label{e-1}
0=(V-\phi)(t_0,i_0)=\max\{(V-\phi)(t,i);(t,i)\in [0,T)\times \S\}.
\end{equation}
Take an arbitrary point $\tilde\mu\in \pb(U)$ and let
\[\mu_t=\tilde \mu, \quad \forall\,t\in [s,T],\]
which is a constant control policy and obviously satisfies the conditions of Definition \ref{def-1}.
According to the dynamic programming principle (Theorem \ref{t2}), we have
$$
V(t_0,i_0)\leq \E\Big[\int_{t_0}^t\!f(r,\La_r,\tilde \mu)\d r+V(t,\La_t)\Big].
$$
Due to \eqref{e-1}, it holds $V\leq \phi$, and hence
\begin{equation}\label{e-2}
\phi(t_0,i_0)\leq \E\Big[\int_{t_0}^t\!f(r,\La_r,\tilde\mu)\d r+\phi(t,\La_t)\Big].
\end{equation}
Applying It\^o-Dynkin's formula to the function $\phi$ (cf. \cite[Theorem 3.1]{GHH}), we get
\begin{equation}\label{e-3}
\E\,\phi(t,\La_t)=\phi(t_0,i_0)+\E\Big[\int_{t_0}^t\!\!\big(\frac{\partial \phi}{\partial r }(r,\La_r)\!+\!Q(\tilde\mu) \phi(r,\La_r) \big)\d r\Big].
\end{equation}
Inserting \eqref{e-3} into \eqref{e-2} leads to
\begin{equation}\label{e-4}
-\E\Big[\int_{t_0}^t\!\!\big(\frac{\partial \phi}{\partial r}(r,\La_r)+Q(\tilde\mu)\phi(r,\La_r)+f(r,\La_r,\tilde\mu )\big)\d r\Big]\leq 0.
\end{equation}
Dividing both sides of \eqref{e-4} by $t-t_0$ and letting $t\downarrow t_0$, we get from the almost sure right-continuity of the trajectories of $(\La_t)$ that
\begin{equation}\label{e-5}
-\frac{\partial \phi}{\partial t}(t_0,i_0)-\sum_{j\neq i_0}q_{i_0j}(\tilde \mu)(\phi(t_0,j)-\phi(t_0,i_0))
+f(t_0,i_0,\tilde \mu)\leq 0.
\end{equation}
Then, by the arbitrariness of $\tilde \mu\in \pb(U)$, $V(t,i)$ is a viscosity subsolution of \eqref{HJB0}.

Next, we proceed to the viscosity supersolution property. Let $(t_0,i_0)\in [0,T)\times \S$ and $\phi\in C^1([0,T)\times \S)$ be a test function such that
\begin{equation}\label{e-6}
0=(V-\phi)(t_0,i_0)=\min\{(V-\phi)(t,i);(t,i)\in [0,T)\times \S\}.
\end{equation}
The desired result will be shown by contradiction. Assume
\begin{equation}\label{e-7}
-\frac{\partial \phi}{\partial t}(t_0,i_0)-\inf_{\mu\in \pb(U)} \big\{Q(\mu)\phi(t_0,i_0)+f(t_0,i_0,\mu)\big\} <0.
\end{equation}
By (H1), the compactness of $\pb(U)$ and the continuity of $f$, we obtain from \eqref{e-7} that there exist  $\veps,\eta>0$ such that for any $0\leq t-t_0\leq \eta$, it holds
\begin{equation}\label{e-8}
-\frac{\partial \phi}{\partial t}(t,i_0)-\inf_{\mu\in \pb(U)} \big\{ Q(\mu) \phi(t,i_0) + f(t,i_0,\mu) \big\} \leq - \veps.
\end{equation}

Let $(t_k)_{k\geq 1}$ be a sequence satisfying $\lim_{k\to \infty} t_k =t_0$.
Using the dynamic programming principle (Theorem \ref{t2}) again, for each $k\geq 1$, there exists   $\alpha^{(k)}=(\La_t^{(k)},\mu_t^{(k)},t_0,i_0)\in \Pi_{t_0,i_0}$  such that
\[ V(t_0,i_0)\geq \E\Big[\int_{t_0}^{\beta_k}\!\!f(r,\La_r^{(k)},\mu_r^{(k)})\d r+V(\beta_k,\La_{\beta_k}^{(k)})\Big]-\frac{\veps}{2}(t_k-t_0),\]
where $\beta_k=t_k\wedge \tau_k$, and $\tau_k$ is defined by
\begin{equation}\label{tau}
\tau_k=\inf\{t\in[t_0,T];\La^{(k)}_t\neq \La^{(k)}_{t_0}\big\} \wedge (t_0+\eta).
\end{equation}
Due to \eqref{e-6}, we have $V\geq \phi$ and
\begin{equation}\label{e-9}
\phi(t_0,i_0)\geq \E\Big[\int_{t_0}^{\beta_k}\!\! f(r,\La_r^{(k)},\mu_r^{(k)})\d r+\phi(\beta_k,\La_{\beta_k}^{(k)})\Big]-\frac{\veps}{2}(t_k-t_0).
\end{equation}
Using It\^o-Dynkin's formula to the function $\phi$, we have
\[\E\Big[\int_{t_0}^{\beta_k}\!\! f(r,\La_r^{(k)},\mu_r^{(k)}) +\big(\frac{\partial \phi}{\partial r} + Q(\mu_r^{(k)})\phi\big)(r,\La_r^{(k)})\d r\Big]\leq \frac{\veps}{2}(t_k-t_0).\]
Then \eqref{e-8} and the definition of $\beta_k$ implies that
\begin{equation}\label{e-10}
\frac{\E[\beta_k-t_0]}{t_k-t_0}\leq \frac{1}{2},\quad k\geq 1.
\end{equation}
On the other hand, by (H1), we have
$$
\p(\beta_k-t_0\leq t_k-t_0) \leq \p\big(\sup_{s\in[t_0,t_k]}|\La_s^{(k)}-\La_{t_0}^{(k)}|>0\big) \leq 1-\e^{-M(t_k-t_0)},
$$
Therefore,
\[\lim_{k\to \infty} \p(\beta_k-t_0\geq t_k-t_0)=1.\]
Since
\[\p(\beta_k-t_0\geq t_k-t_0)\leq \frac{\E[\beta_k-t_0]}{t_k-t_0}\leq 1,\]
we get finally that
\begin{equation}\label{e-11}
\lim_{k\to\infty}\frac{\E[\beta_k-t_0]}{t_k-t_0}=1,
\end{equation}
which contradicts \eqref{e-10}. Consequently, we have
\begin{equation}\label{e-12}
-\frac{\partial \phi}{\partial t}(t_0,i_0)-\inf_{\mu\in \pb(U)}
\Big\{\sum_{j\neq i_0}\big(\phi(t_0,j)-\phi(t_0,i_0)\big)+f(t_0,i_0,\mu)\Big\}\geq 0.
\end{equation}
This means that $V(t,i)$ is a viscosity supersolution of \eqref{HJB0}. We conclude the proof of this theorem by the definition of viscosity solution to \eqref{HJB0}.
\end{proof}

In the end let us discuss the uniqueness of the viscosity solution to  \eqref{HJB0}.
For this purpose it is sufficient to establish the following comparison principle  for \eqref{HJB0}. We shall develop the method used to establish the comparison principle for HJB equations associated with diffusion processes to the equations associated with purely jumping processes.

\begin{theorem}\label{com-t}
Assume  the conditions of Proposition \ref{prop1} hold. Let $V_1$ (resp. $V_2$) be a viscosity supersolution (resp. viscosity subsolution) of \eqref{HJB0} in $[0, T) \times \S$. Then
$$
\sup_{[0, T] \times \S} [V_2 - V_1] = \sup_{\{T\} \times \S} [V_2 - V_1]=0.
$$
\end{theorem}

\begin{proof} Obviously, we just need to show that
\begin{equation}\label{supp0}
\sup_{[0, T] \times \S} [V_2 - V_1] \leq \sup_{\{T\} \times \S} [V_2 - V_1]=0.
\end{equation}
Under the assumptions of Proposition \ref{prop1}, the boundedness of $f$ and $g$ implies the boundedness of the value function. Therefore, there exists a constant $K_0>0$ such that
\begin{equation}
K_{0}\geq \sup_{t \in [0, T]}\sup_{j\in\S}\big\{ |V_1 (t, j)| \vee |V_2 (t, j)|\big\}.
\end{equation}
There exists a sequence of $C^2(\R)$ functions $\lambda_n(x)$ such that $\lambda_n(x)=0$ for $x\leq 0$, $0<\lambda_n'(x)<1$, $\lambda_n(x)\uparrow \max\{x,0\}$ as $n\to \infty$.
Let
\[\eta_n(s,t)=t+\lambda_n(s\!-\!t),\quad s,t\in [0,T].\]
Then $\eta_n(s,t)\uparrow \max\{s,t\}$ as $n\to \infty$.
Define a function on $[0, T] \times [0, T]$ as
$$
\Psi_{i_0}^{n} (t, s) = V_2 (t, i_0) - V_1 (s, i_0) - \frac{1}{2 \delta} (t- s)^2 + \frac{\beta}{\delta} (\eta_n(s,t) - T),
$$
where $\delta, \beta > 0$ are two parameters. Again, the continuity of $V_1$ and $V_2$ implies that $\Psi_{i_0}^n$ achieves the maximum on $[0, T] \times [0, T]$. Denote  by $(\bar{t}, \bar{s})\in [0,T]\times [0,T]$ an arbitrary one of the maximum points, and note that $(\bar t,\bar s)$ may depend on the parameters $\delta,\beta$.

We first give an estimate of the distance between $\bar s$ and $\bar t$.  For any $\rho \geq 0$, let
\begin{align*}
D_{\rho} &= \left\{(t, s) \in [0, T] \times [0, T] : |t-s|^2 \leq \rho \right\}, \\
m_{i_0}^{(1)} (\rho) &= 2 \sup \left\{ |V_1 (t, i_0) - V_1 (s, i_0)| : (t, s) \in D_{\rho} \right\}, \\
m_{i_0}^{(2)} (\rho) &= 2 \sup \left\{ |V_2 (t, i_0) - V_2 (s, i_0)| : (t, s) \in D_{\rho} \right\}.
\end{align*}
Then $m_{i_0}^{(1)}$ and $m_{i_0}^{(2)}$ are increasing functions satisfying $m_{i_0}^{(1)} (0) = m_{i_0}^{(2)} (0) =0$. Moreover, it follows from the continuity of $V_1$ and $V_2$ and the compactness of $[0, T] \times [0, T]$  that $m_{i_0}^{(1)}, m_{i_0}^{(2)} $ are continuous. Since $V_1 (\cdot, i_0)$ and $V_2 (\cdot, i_0)$ are bounded, $m_{i_0}^{(1)}$ and $m_{i_0}^{(2)}$ are  bounded as well and  denoted by $M_{i_0} := \sup \{m_{i_0}^{(1)} (\rho) + m_{i_0}^{(2)} (\rho) : \rho \geq 0\} < \infty$. We obtain from the fact $\Psi_{i_0}^n (\bar{t} \vee \bar{s}, \bar{t} \vee \bar{s}) \leq \Psi_{i_0}^n (\bar{t}, \bar{s})$ that
$$
\frac{1}{\delta} (\bar{t} - \bar{s})^2 \leq 2 \left( V_2 (\bar{t}, i_0) - V_2 (\bar{t} \vee \bar{s}, i_0) + V_1 (\bar{t} \vee \bar{s}, i_0) - V_1 (\bar{s}, i_0) \right) \leq M_{i_0}.
$$
Hence,
\begin{equation}\label{supp1}
|\bar{t} - \bar{s}| \leq \sqrt{\delta M_{i_0}}, \quad \text{and hence}\ \ \bar t-\bar s\to 0,\ \ \text{as}\ \delta\to 0.
\end{equation}

Next, we shall show by contradiction that $\bar t$ equals to $T$. Assume that  $\bar{t}\in [0, T)$. Define an auxiliary function on $[0, T] \times \S$ as
$$
\psi_{i_0}^{(1)} (s, j) = - \frac{1}{2 \delta} (\bar{t} - s)^2 - 2K_{0}  \big(1- \mathbf{1}_{i_0} (j)\big) + \frac{\beta}{\delta} (\eta_n(s,\bar t) - T).
$$
For each $s \in [0, T]$, since $\Psi^n_{i_0} (\bar{t}, s) \leq \Psi^n_{i_0} (\bar{t}, \bar{s})$, it holds that
\begin{equation*}
V_1 (\bar{s}, i_0) + \frac{1}{2 \delta} (\bar{t} - \bar{s})^2 - \frac{\beta}{\delta} (\eta_n(\bar s,\bar t) - T) \leq V_1 (s, i_0) + \frac{1}{2 \delta} (\bar{t} - s)^2 - \frac{\beta}{\delta} (\eta_n(s,\bar t) - T),
\end{equation*}
and further for each $j \in \S$, $j \neq i_0$,
\begin{equation*}
\begin{split}
2K_{0} &\geq V_1(s,i_0)-V_1(s,j)\\
& \geq V_1 (\bar{s}, i_0) - V_1 (s, j) + \frac{1}{2 \delta} (\bar{t} - \bar{s})^2 - \frac{1}{2 \delta} (\bar{t} - s)^2 - \frac{\beta}{\delta} (\eta_n(\bar s,\bar t)-\eta_n(s,\bar t)).
\end{split}
\end{equation*}
Hence, $(\bar{s}, i_0)$ is the minimum point of the function $(s,j)\mapsto V_1(s,j)-\psi_{i_0}^{(1)}(s,j)$. Since $V_1$ is the viscosity supersolution of \eqref{HJB0}, we have
\begin{equation}\label{supp4}
- \frac{1}{\delta} (\bar{t}-\bar{s}) - \frac{\beta}{\delta } \lambda_n'(\bar s-\bar t)   - \inf_{\mu \in \pb(U)} \left\{ - 2K_{0}  q_{i_0} (\mu) + f (\bar{s}, i_0, \mu)\right\} \geq 0.
\end{equation}
Similarly, consider the test function on $[0,T] \times \S$ as
$$
\psi_{i_0}^{(2)} (t, j) = \frac{1}{2 \delta} (t - \bar{s})^2 + 2K_{0}  \big(1- \mathbf{1}_{i_0} (j)\big) - \frac{\beta}{\delta} (\eta_n(\bar s,t) -T).
$$
Then, $\Psi_{i_0}^n(t,\bar s)\leq \Psi_{i_0}^n(\bar t,\bar s)$ implies that for each $t \in [0, T]$,
\begin{equation*}
V_2 (t, i_0) - \frac{1}{2 \delta} (t - \bar{s})^2 + \frac{\beta}{\delta } (\eta_n(\bar s, t) - T) \leq V_2 (\bar{t}, i_0) - \frac{1}{2 \delta} (\bar{t} - \bar{s})^2 + \frac{\beta}{\delta} (\eta_n(\bar s,\bar t) - T) ,
\end{equation*}
and for each $j \in \S$ with $j \neq i_0$,
\begin{equation*}
2K_{0}  \geq V_2 (t, j) - V_2 (\bar{t}, i_0) + \frac{1}{2 \delta} (\bar{t} - \bar{s})^2 - \frac{1}{2 \delta} (t - \bar{s})^2 + \frac{\beta}{\delta} ( \eta_n(\bar s,t)-\eta_n(\bar s,\bar t)).
\end{equation*}
This  means that $(\bar{t}, i_0)$ is a maximum point of $(t,j)\mapsto V_2(t,j) - \psi_{i_0}^{(2)}(t,j)$. Since $V_2$ is the viscosity subsolution of \eqref{HJB0}, we have
\begin{equation}\label{supp5}
 \frac{\beta}{\delta}\big(1-\lambda_n'(\bar s-\bar t)\big) - \frac{1}{\delta} (\bar{t}-\bar{s}) - \inf_{\mu \in \pb(U)} \left\{ 2K_0 q_{i_0} (\mu) + f (\bar{t}, i_0 , \mu)\right\} \leq 0.
\end{equation}
Combining the inequalities (\ref{supp1}), (\ref{supp4}), (\ref{supp5}) and \eqref{d-0}, we arrive at
\begin{equation}\label{supp9}
\begin{split}
  \frac{\beta}{\delta}&\leq \inf_{\mu\in\pb(U)} \big\{2K_0q_{i_0}(\mu)+f(\bar t,i_0,\mu)\big\}-\inf_{\mu\in\pb(U)}\big\{-2K_0q_{i_0}(\mu)+f(\bar s,i_0,\mu)\big\}\\
  &=\sup_{\mu\in\pb(U)}\big\{ 2K_0q_{i_0}(\mu)-f(\bar s,i_0,\mu)\big\}-\sup_{\mu\in\pb(U)}\big\{-2K_0q_{i_0}(\mu)-f(\bar t,i_0,\mu)\big\}\\
  &\leq \sup_{\mu\in\pb(U)}\big\{4K_0q_{i_0}(\mu)+f(\bar t,i_0,\mu)-f(\bar s,i_0,\mu)\big\}\\
  &\leq 4K_0 M+C_0|\bar t-\bar s|.
\end{split}
\end{equation}
Invoking the estimate \eqref{supp1}, this yields that
\[\beta\leq 4K_0M\delta+C_0\delta^{3/2}\sqrt{M_{i_0}}.\]
Thus, letting $\delta\to 0$, we get that $\beta\leq 0$, which contradicts the assumption that $\beta>0$.
Consequently, it must hold
\begin{equation}
\bar t=T .
\end{equation}

By the choice of $(\bar t,\bar s)$, it holds that  for every $t\in [0,T)$,
\begin{equation}\label{sup-1}
\begin{aligned}V_2(t,i_0)-V_1(t,i_0)+\frac{\beta}{\delta}(t-T)&
=\Psi_{i_0}^n(t,t)\leq \Psi_{i_0}^n(\bar t,\bar s)\\
&=V_2(T,i_0)-V_1(\bar s,i_0)-\frac{1}{2\delta}(T-\bar s)^2\\
&\leq V_2(T, i_0)-V_1(\bar s,i_0).
\end{aligned}
\end{equation}
Thus, letting first $\beta\to 0$ and then $\delta\to 0$, noting $\lim_{\delta\to 0}|\bar t-\bar s|=0$ due to \eqref{supp1}, we obtain that
\[V_2(t,i_0)-V_1(t,i_0)\leq V_2(T,i_0)-V_1(T, i_0).\]
The desired conclusion \eqref{supp0} follows from the arbitrariness of $i_0\in \S$.
\end{proof}


The following uniqueness result is an immediate result of Theorem \ref{t3} and Theorem \ref{com-t}.
\begin{corollary}\label{t4}
Under the conditions of Proposition \ref{prop1}, the value function $V(t,i)$ is the unique viscosity solution to the equation \eqref{HJB0}.
\end{corollary}

Next, noticing that the equation \eqref{HJB0} does not rely on the delay-dependent control policies, we shall take advantage of this property to show the existence of an optimal Markovian control policy over the class of delay-dependent controls.

\begin{theorem}\label{t5}
Under the conditions of Proposition \ref{prop1}, for every $t\in [0,T],\, i\in\S$, there exists an optimal control $\alpha^\ast$  for $V(t,i)$, which depends only on the current state of the process $(\La_t)$, i.e. a Markovian control policy.
\end{theorem}

\begin{proof}
  Introduce a sub-class $\Pi_{s,i}^m$ of $\Pi_{s,i}$ by
  \[\Pi_{s,i}^m=\big\{\alpha\in (\La_t,\mu_t,s,i)\in \Pi_{s,i}; \ \exists\,h:\S\to \pb(U)\ \text{such that $\mu_t=h(\La_t)$}\big\},\]
  which is the class of stationary randomized Markov  policy. Let
  \begin{equation}\label{e-cor-1}
  \wt V(s,i)=\inf_{\alpha\in \Pi_{s,i}^m}J(s,i,\alpha),
  \end{equation} which is consistent with the value function studied in \cite[p.1069]{GHH}.
  According to \cite[Theorem 4.1]{GHH} and Proposition \ref{prop1}, for each $i \in \S$, $s\mapsto \widetilde{V} (s, i)$ is differentiable almost everywhere for every $i\in\S$ and satisfies the HJB equation \eqref{HJB0} almost everywhere. Moreover, $\widetilde{V} (s, i)$ is also a viscosity solution to \eqref{HJB0}. Hence, the uniqueness of viscosity solution given in Theorem \ref{t4} means that $\widetilde{V}(s, i) = V(s, i)$.
  Using  \cite[Theorem 4.1]{GHH} again or along the procedure of Theorem \ref{t1}, there exists an $\tilde \alpha\in \Pi_{s,i}^m$ such that $\wt V(s,i)=J(s,i,\tilde \alpha)$.  Therefore,
  \begin{equation}\label{e-cor-2}
  V(s,i)=\wt V(s,i)=J(s,i,\tilde \alpha),
  \end{equation}
  which means that $\tilde \alpha\in \Pi_{s,i}^m\subset \Pi_{s,i}$ is  the desired optimal control policy in $\Pi_{s,i}$ associated with $V(s,i)$.
\end{proof}

\section*{Appendix}
\label{Appendix}

In the section, we construct an example to illustrate that for discrete-time decision processes in an infinite state space, the optimization problem may have essential difference between the control mechanism over history-dependent control policies and over Markovian policies. In this example, the value function corresponding to taking infimum over Markovian policies equals $+\infty$, while the one over history-dependent policies equals $-\infty$. Therefore, when analyzing the influence of control policy class on the value function, more attentions should be paid.

Let the state space $X=\Z_+=\{0, 1,2,\ldots\}$ and the action space $A=\Z_+$. Denote $\pb(A)$ the set of probability measures on $A$.  Consider the transition probability matrices given by
\begin{gather*}
  P_0(0)=1,\qquad P_1(j|i,a)=\begin{cases}
                               \frac{1}{K j^2}, & \mbox{$j\neq 0$,}   \\
                              0, & \mbox{$j=0$}.
                             \end{cases} \\
                             P_2(0|i,a)=P_3(0|i,a)=1,\qquad \forall\, i\in X,\ a\in A.
\end{gather*} Here $K=\sum_{j=1}^\infty \frac1{j^2}$ is a constant.
Let $\{\xi_k;k=0,1,2,3\}$ denote the controlled process. By the definition of $P_t(\cdot|i,a)$ above,  it holds that
\[\p(\xi_0=0)=1, \ \p(\xi_1\geq 1)=1,\ \p(\xi_2=0)=\p(\xi_3=0)=1.\]

For a probability measure $\mu$ on $A$, denote by
\[m_1(\mu)=\sum_{i\geq 0} i\mu(i),\quad m_2(\mu)=\sum_{i\geq 0} i^2\mu(i),\quad \mathrm{var}(\mu)=m_2(\mu)-(m_1(\mu))^2.\]
Let
\begin{equation}\label{ee-1} \rho(x,\mu)=\big(-2m_1(\mu)+\infty\cdot \mathrm{var}(\mu)\big)\mathbf{1}_{\{m_1(\mu)<+\infty\}}
\end{equation}
 for $x\in X$ and $\mu\in \pb(A)$ and using the convention $\infty \cdot 0=0$. By \eqref{ee-1}, $\rho(x,\mu)$ takes value in $(-\infty,+\infty]$, and if $m_1(\mu)<+\infty$ and $\mathrm{var}(\mu)>0$, then $\rho(x,\mu)=+\infty$. So, if $\rho(x,\mu)<+\infty$, $\mu$ must be a Dirac measure in the form $\mu=\delta_{i_0}$ for some $i_0\in A$, which plays an important role below.

Define the cost function $c_t(\cdot,\cdot)$ by
\begin{equation}\label{ee-2}
\text{$c_1(i,\mu)=0$, $c_2(i,\mu)=i$, $c_3(i,\mu)=\rho(i,\mu)$ \ for $i\in X$ and $\mu\in\pb(A)$}.
\end{equation}
For the control policy $\pi$,
\begin{equation}\label{ee-3}
V^\pi(i):=\E_i\Big[ \sum_{t=1}^3 c_t(\xi_{t-1},\pi_t)\Big].
\end{equation}
Let $\Pi$ be the set of all history-dependent control policies, and $\Pi^M$ the set of all Markov control policies. Clearly,  $\Pi^M\subset \Pi$. The corresponding value functions are given by
\[V(i)=\inf_{\pi\in \Pi} V^{\pi}(i), \quad V^M(i)=\inf_{\pi\in \Pi^M} V^\pi(i),\quad i\in X.\]
We shall show that
\begin{equation}\label{ee-4}V(0)=-\infty,\quad\ \text{but}\ \ V^M(0)=+\infty.\end{equation}
Indeed, according to \eqref{ee-2},
 \begin{align*}
   V^\pi(0)&=\E\big[c_2(\xi_1,\pi_2)+c_3(\xi_2,\pi_3)\big]=\E[\xi_1+\rho(\xi_2,\pi_3)\big]\\
   &=\E\big[\xi_1+\rho(0,\pi_3)\big]\qquad \qquad (\text{as $\xi_2=0$ a.s.}).
 \end{align*}
Note that $\E[\xi_1]=\sum_{i=1}^\infty \frac{1}{K i}=+\infty$.

For every Markov control policy $\pi$,
 \begin{itemize}
   \item[(i)] If $m_1(\pi_3)<\infty$, and $\mathrm{var}(\pi_3)>0$, then $\rho(0,\pi_3)=+\infty$, and hence $V^\pi(0)=+\infty$.
   \item[(ii)] If $\rho(0,\pi_3)<+\infty$, then $\pi_3$ must be a Dirac measure. In addition, since $\xi_2\equiv 0$ and $\pi_3$ is a Markov control policy, which is a functional of $\xi_2$, there exists a function $f:X\to A$ such that
       \[\pi_3(\d x)=\delta_{f(0)}(\d x).\]
       Hence,
        $\dis V^\pi(0)=\E[\xi_1-2f(0)]=+\infty.$
 \end{itemize}
From the discussion in (i) and (ii) above, we obtain that \[V^M(0)=\inf_{\pi\in \Pi^M}V^\pi(0)=+\infty.\]

For the set of history-dependent control policies, we choose a special one $\tilde \pi$ given by
\[\tilde \pi_1(\d x)=\delta_0(\d x),\quad \tilde \pi_2(\d x)=\delta_0(\d x),\ \ \tilde \pi_3(\d x)=\delta_{\xi_1}(\d x).\] Note that $\tilde \pi_3$ depends on $\xi_1$, not on $\xi_2$, so $\tilde \pi$ is not a Markov control policy.
Then
\[V^{\tilde\pi}(0)=\E[\xi_1-2\xi_1]=-\E[\xi_1]=-\infty.\]
Hence,
\[V(0)=\inf_{\pi\in \Pi}V^\pi(0)\leq V^{\tilde\pi}(0)=-\infty.\]
Consequently, we have proved the desired result \eqref{ee-4}.

\section*{Acknowledgements}

This work was supported in part by National Key R\&D Program of China (No. 2022YFA1000033), the National Natural Science Foundation of China (No.~12271397, 11831014) and the Guangdong Basic and Applied Basic Research Foundation (No.~2022A1515010222).

\section*{Declarations}

\noindent {\bf Conflict of Interests:} The authors declare that they have no conflict of interest.





\end{document}